\documentclass[11pt, reqno]{amsart}
\usepackage{amsmath, amsthm, amscd, amsfonts, amssymb, graphicx, color}
\usepackage[bookmarksnumbered,plainpages]{hyperref}

\newtheorem{thm}{Theorem}[section]
\newtheorem{cor}[thm]{Corollary}
\newtheorem{lem}[thm]{Lemma}

\theoremstyle{definition}
\newtheorem{defn}[thm]{Definition}
\theoremstyle{remark}
\newtheorem{rem}[thm]{Remark}
\numberwithin{equation}{section}

\newcommand{\A}{\mathcal{A}}
\newcommand{\f}{\mathcal{F}}
\newcommand{\LL}{\mathcal{L}}

\newcommand{\B}{\mathcal{B}}

\begin{document}
\title[]{Filters and Semigroup compactification properties }
\author{M. Akbari Tootkaboni}
\address{Department of Mathematics, Faculty of science, Shahed
University of Tehran\\
P. O. Box:   18151-159, Tehran, Iran}
\email{akbari@shahed.ac.ir}
\author{H.\ R.\ E.\ Vishki}
\address{Department of Pure Mathematics and Centre of Excellence
in Analysis on Algebraic Structures (CEAAS), Ferdowsi University of Mashhad\\
P. O. Box 1159, Mashhad 91775, Iran.} \email{vishki@ferdowsi.um.ac.ir}

\subjclass[2000]{}
\keywords{}
\maketitle

\begin{abstract}
Stone-$\breve{\mbox{C}}$ech compactifications derived from a
discrete semigroup $S$ can be considered as the spectrum of the
algebra $\mathcal{B}(S)$ or as a  collection of ultrafilters on $S$.
What is certain and indisputable is the fact that filters play an
important role in the study of Stone-$\breve{\mbox{C}}$ech
compactifications
derived from a discrete semigroup.\\
It seems that filters can play a role in the study of general
semigroup compactifications too. In the present paper, first
we review the characterizations of semigroup compactifications
in terms  of
filters and then extend some of the results in
\cite{Talin1} concerning the
Stone-$\breve{\mbox{C}}$ech compactification  to
a semigroup compactification associated with  a Hausdorff
semitopological semigroup.
\end{abstract}
\maketitle
\section{Introduction}
It is well known that for a discrete semigroup $(S,.)$, the
Stone-$\breve{\mbox{C}}$ech compactification $\beta S$ can be
regarded as the set of all ultrafilters on $S$ and so the semigroup
operation of $S$ can be extended in a natural way to $(\beta S,+)$.
Indeed, for $p,\,q\in \beta S$, $p+q\in\beta S$ is defined by
$p+q=\{A\subseteq S:\Omega_{q}(A)\in p\}$, where
$\Omega_{q}(A)=\{x\in S:\lambda_{x}^{-1}(A)\in q\}$ and
$\lambda_{x}$ denotes the left translation by $x$. Moreover
$\{\overline{A}:A\subseteq S\}$ forms a basis for the topology of
$\beta S$, where $\overline{A}=\{p\in \beta S:A\in p\}$, and $(\beta
S,+)$ is a compact Hausdorff right topological semigroup with $S$
contained in its topological center. The topological and algebraic
structures of $\beta S$ from the ultrafilters point of view have
been studied extensively in the book of Hindman and Strauss
\cite{hindbook}. For an extensive account about ultrafilters the
readers may refer to \cite{Ultra} and \cite{Gil}. Alternatively $\beta S$
viewed as a semigroup compactification of $S$ may be obtained
as the set of all multiplicative means on $\mathcal{B}(S)$
(see \cite{Analyson}).

 For an m-admissible subalgebra $\f$ (of
$\mathcal{CB}(S)$) on a Hausdorff semitopological semigroup $S$, it
has been shown in \cite{Akbari}, that the semigroup compactification
$S^{\f}$(which is the set of all multiplicative means on $\f$) can
be obtained  as a suitable quotient space of the set of all
$z$-ultrafilters.

The purpose of this paper is to use filter methods mentioned in
\cite{Talin1} to obtain more information about $S^{\f}$. Here, for
filters obtained from an m- admissible subalgebra $\f$, the
operation of addition for filters has been defined.  In the special
case in which $S$ is a discrete semigroup the definition will be the
same as the definition ascribed to Glazer (see
\cite{hindbook}). Also we give some characterizations of closed
subsemigroups and closed (left and right) ideals of $S^{\f}$, whose
discrete versions are given by Papazian in \cite{Talin1}. In fact,
in this paper Lemma 3.9, Theorems 3.13, 3.15,  3.18, 3.20 and Corollaries 3.14, 3.16   3.19 have
been stated by Papazian in \cite{Talin1}, when $S$ is a discrete
semigroup.

\section{Preliminaries}

 Throughout this paper $S$ is a Hausdorff semitopological semigroup
 and $\mathcal{CB}(S)$
 is the $\mathcal{C}^{*}$-algebra of all  bounded continuous
 complex-valued functions on
 $S$ with supremum norm. Let $\mathcal{F}$ be an m-admissible subalgebra
 of $\mathcal{CB}(S)$, then the set of all multiplicative means of
 $\f$, which we denote by $S^{\f}$, equipped with the
 Gelfand topology is a compact Hausdorff right topological semigroup,
 whose multiplication is given by $\mu\nu(f)=\mu( \mathbf{T}_{\nu}f)$, where
$\mathbf{T}_{\nu}(f)(s)=\nu(L_{s}f)$ for all $s\in S$, $\mu,\nu\in
S^{\mathcal{F}}$ and $f\in\f$. Furthermore the evaluation mapping
$\varepsilon :S\rightarrow S^{\mathcal{F}}$ is given by
$\varepsilon(s)(f)=f(s)$ and  is a continuous homomorphism onto a dense
subsemigroup of $S^{\mathcal{F}}$ contained in the topological
center of $S^{\f}$. In other words, $(\varepsilon, S^{\f})$ is a
semigroup compactification of $S$. Also
$\varepsilon^{*}:\mathcal{C}(S^{\mathcal{F}})\rightarrow\mathcal{F}$
is an isometric isomorphism and
$\widehat{f}=(\varepsilon^{\ast})^{-1}(f)\in
\mathcal{C}(S^{\mathcal{F}})$ for $f\in\f$ is given by
$\widehat{f}(\mu)=\mu(f)$ for all $\mu\in S^{\f}$, (see
\cite{Analyson}).

 Now we quote some prerequisite material from
\cite{Akbari} for the description of $S^{\f}$ in terms of filters.
 For all $f\in\f$, $Z(f)=f^{-1}(\{0\})$ is called a zero set and
 we denote the collection of all zero sets by $Z(\f)$.
\begin{defn}
$\mathcal{A}\subseteq Z(\mathcal{F})$ is called a $z$-filter if,\\
(i) $\emptyset\notin \mathcal{A}$ and $S\in\A,$\\
(ii) if $A,B\in \mathcal{A}$, then $A\bigcap B\in\mathcal{A},$\\
(iii) if $A\in \mathcal{A}$, $B\in Z(\mathcal{F})$ and $A\subseteq
B$ then $B\in \mathcal{A}$.
\end{defn}
 A $z$- filter is said to be a $z$-ultrafilter if it is not contained
properly in any other $z$-filter. We put $\mathcal{F}S=\{p\subseteq Z(\f): p\,\,\mbox{is an
$z$-ultrafilter}\}$. It is obvious that
$\widehat{x}=\{Z(f):f\in\f,\,\,f(x)=0\}$ is an ultrafilter and
$\widehat{x}\in\f S$. If $p\in\f S$, then there is $\mu\in
S^{\mathcal{F}}$ such that $\bigcap_{A\in
p}\overline{\varepsilon(A)}=\{\mu\},$ also for each $\mu$ in
$S^{\mathcal{F}}$ there exists $p\in \f S$ such that $\bigcap_{A\in
p}\overline{\varepsilon(A)}=\{\mu\}.$  If $A\in Z( \mathcal{F})$,
and there exists a neighborhood $U$ of $\mu$ in $S^{\mathcal{F}}$
such that $U\subseteq \overline{\varepsilon(A)}$, then $A\in p$
(see \cite{Akbari}).

 Unlike the discrete
setting, there is no simple relationship between $S^{\f}$ and $\f
S$. $\f S$ is equipped with a topology whose base is $\{
(\widehat{A})^{c}: A\in Z(\f)\}$, where $\widehat{A}=\{p\in\f S:
A\in p\}$ is a compact space which is not Hausdorff in general.
\begin{defn}
We define the relation $\sim$ on $\f S$ by putting $p\sim q$ if
$$\bigcap_{A\in p}\overline{\varepsilon(A)}=\bigcap_{B\in
q}\overline{\varepsilon(B)}.$$
\end{defn}
It is obvious that $\sim$ is an equivalence relation on $\f S$. Let
$[p]$ be the equivalence class of $p\in \f S$, and let $\frac{\f
S}{\sim}$ be the corresponding quotient space with the quotient map
$\pi:\f S\rightarrow \frac{\f S}{\sim}$. For every $p\in \f S$
define $\widetilde{p}$ by $\widetilde{p}=\bigcap [p]$, put
$\widetilde{A}=\{\widetilde{p}:A\in p\}$ for $A\in Z(\f)$ and
$\mathcal{R}=\{\widetilde{p}:p\in\f S\}$. It is obvious that
 $\{(\widetilde{A})^{c}:A\in Z( \mathcal{F})\}$ is a basis for a
topology on $ \mathcal{R}$, $\mathcal{R}$ is a Hausdorff and compact
space and also $S^{\mathcal{F}}$ and $\mathcal{R}$  are homeomorphic
(see \cite{Akbari}). So we have
$\mathcal{R}=\{\mathcal{A}^{\mu}:\mu\in S^{\mathcal{F}}\}$, where
$\mathcal{A}^{\mu}=\widetilde{p}$.
\begin{defn}
For all $x,y\in S$, we define
$$\mathcal{A}^{\varepsilon(x)}*\mathcal{A}^{\epsilon(y)}=\{ Z(f)\in Z(\mathcal{F}):\;
Z(\mathbf{T}_{\varepsilon(y)}f)\in \mathcal{A}^{\varepsilon(x)}\}.$$
\end{defn}
\begin{lem}
$\mathcal{A}^{\varepsilon(x)}*\mathcal{A}^{\varepsilon(y)}=\mathcal{A}^{\varepsilon(xy)}$,
for each $x,y\in S$
\end{lem}
{\bf Proof :} If $Z(f)\in \mathcal{A}^{\varepsilon(xy)}$, then
$\varepsilon(xy)\in\overline{\varepsilon(Z(f))}$, and so $f(xy)=0$
and this implies that $\mathbf{T}_{\varepsilon(y)}f(x)=f(xy)=0$.
Therefore $x\in Z(\mathbf{T}_{\varepsilon(y)}f)$ and
$Z(\mathbf{T}_{\varepsilon(y)}f)\in \mathcal{A}^{\varepsilon(x)}$ so
$Z(f)\in \mathcal{A}^{\varepsilon(x)}*\mathcal{A}^{\varepsilon(y)}$.
Also if
$Z(f)\in\mathcal{A}^{\varepsilon(x)}*\mathcal{A}^{\varepsilon(y)}$
then $Z(\mathbf{T}_{\varepsilon(y)}f)\in
\mathcal{A}^{\varepsilon(x)}$, and so
$\mathbf{T}_{\varepsilon(y)}f(x)=f(xy)=0$. This implies that
$Z(f)\in \mathcal{A}^{\varepsilon(xy)}$. Therefore
$\mathcal{A}^{\varepsilon(x)}*\mathcal{A}^{\varepsilon(y)}=\mathcal{A}^{\varepsilon(xy)}$.$\hfill\Box$
\begin{defn}
Let $\{x_{\alpha}\}$ and  $\{y_{\beta}\}$ be two nets in $S$, such
that $lim_{\alpha}\varepsilon(x_{\alpha})=\mu$ and
$lim_{\beta}\varepsilon(y_{\beta})=\nu$, for $\mu,\nu\in
S^{\mathcal{F}}$. We define
$$\mathcal{A}^{\mu}*\mathcal{A}^{\nu}=lim_{\alpha}(lim_{\beta}(
\mathcal{A}^{\varepsilon(x_{\alpha})}*\mathcal{A}^{\varepsilon(y_{\beta})})).$$
 \end{defn}
 It is obvious that Definition 2.5 is well-defined and $(\mathcal{R},e)$ is
 a compact right topological semigroup,
where $e:S\rightarrow \mathcal{R}$ is defined by
$e[x]=\widehat{x}$. Also the mapping
$\varphi:S^{\mathcal{F}}\rightarrow \mathcal{R}$ defined by
$\varphi(\mu)=\widetilde{p}$, where $\bigcap_{A\in
p}\overline{\varepsilon(A)}=\{\mu\}$, is an isomorphism (see
\cite{Akbari}).

The operation ``$\cdot$" on S extends uniquely to $(\mathcal{R},*)$. Thus
$(\mathcal{R}, e)$ is a semigroup compactification of $(S,\cdot)$, that
$e:S\rightarrow \mathcal{R}$, and $e[x]=\{A\in Z(\f)\;:x\in A\}$
for each $x\in S$, is an evaluation map. Also $e[S]$ is a subset of the
topological center of $\mathcal{R}$ and
$cl_{\mathcal{R}}(e[S])=\mathcal{R}$. For more details see
\cite{Akbari}.

Hence $S^{\f}$ and $\mathcal{R}$ are topologically isomorphic and
so $S^{\f}\simeq\mathcal{R}$.

 In this paper we write $\overline{A}$ for
$\widetilde{A}$, and we define
$$[S]^{<\omega}=\{A\subseteq S: \mbox{$A$ is a finite and nonempty
subset of $S$}\}.$$ Also $\lambda_{x}$ denotes the left translation
by $x\in S$ and $\rho_{\mu}$  denotes the right translation by
$\mu\in S^{\f}$.
\begin{lem}
Let $A,\,\,B\in Z(\f)$. Then $(\overline{A})^{\circ}\cap
(\overline{B})^{\circ}=(\overline{A\cap B})^{\circ}$.
\end{lem}
{\bf Proof :} $(\overline{A\cap
B})^{\circ}\subseteq(\overline{A})^{\circ}\cap
(\overline{B})^{\circ}$ follows from basic topology.

 For the
converse, let $x\in(\overline{A})^{\circ}\cap
(\overline{B})^{\circ}$, then there exists $f\in\f$ such that $x\in
(\overline{Z(f)})^{\circ}\subseteq \overline{Z(f)}\subseteq
(\overline{A})^{\circ}\cap (\overline{B})^{\circ}$. Hence
$Z(f)\subseteq A\cap B$ and so $x\in
(\overline{Z(f)})^{\circ}\subseteq \overline{Z(f)}\subseteq
(\overline{A\cap B})$ and this implies that $x\in (\overline{A\cap
B})^{\circ}$.$\hfill\Box$\\
{\bf Remark :} If $A,\,\,B\in Z(\f)$, then necessarily
$\overline{A}\cap \overline{B}=\overline{A\cap B}$ is not true. For
example let $\f$ be an m- admissible subalgebra of
$\mathcal{CB}(\mathbb{R})$ such that $S^{\f}=\mathbb{R}_{\infty}$.
Then there exist $A,\,B\in Z(\f)$ such that $A\cap B=\emptyset$,
$\infty\in\overline{A}$ and $\infty\in\overline{B}$.

\section{Some results on semigroup compactifications}

When $S$ be a discrete semigroup, operations can be defined on the
collection of filters on $S$ which are helpful in the study of
$\beta S$. When $S$ is a Hausdorff semitopological semigroup and
$\f$ is an m-admissible subalgebra of $\mathcal{CB}(S)$, we shall
consider $S^{\f}$ as a collection of filters. Finding the rules
governing them is the essential goal of this paper.
 \begin{defn}
 A filter $\mathcal{A}\subseteq Z(\mathcal{F})$ is called a pure filter if
 for some $p\in\f S$, $\mathcal{A}\subseteq p$
implies that $\mathcal{A}\subseteq \widetilde{p}$. The collection of
all pure filters is denoted by $\mathcal{P}(\f)$.
\end{defn}
It is obvious that $\widetilde{p}\in\mathcal{R}$ is a maximal member
of $\mathcal{P}(\f)$. We want to know which semigroup properties are
suggested by filters, so the subcollection of all pure filters can
be important.
\begin{defn}
For a filter $\mathcal{A}\subseteq Z(\mathcal{F})$, we define
$$\overline{\mathcal{A}}=\{\widetilde{p}\in \mathcal{R}:\mbox{There
exists $p\in\f S$ such that }\mathcal{A}\subseteq p\}.$$
\end{defn}
 \begin{lem}
 Let $\A$ and $\B$ be filters. Then the following statements hold.\\
 (i) $\overline{\A}=\bigcap_{A\in
 \Gamma}(\overline{A})^{\circ}=\bigcap_{A\in
 \Gamma}\overline{A}=\bigcap_{A\in\A}\overline{A}$,
  where $\Gamma=\{A\in Z(\f): \overline{\A}\subseteq
 (\overline{A})^{\circ}\}$.\\
(ii) For a filter $\mathcal{A}\subseteq Z(\f)$,
$\overline{\mathcal{A}}$ is a closed subset of $\mathcal{R}$.\\
(iii) $\A\subseteq\B$ then $\overline{\B}\subseteq\overline{\A}$.\\
(iv) If $\mathcal{J}$ is a closed subset of $\mathcal{R}$. Then
$\mathcal{A}=\bigcap \mathcal{J}$ is a pure filter and
$\overline{\mathcal{A}}=\mathcal{J}$.\\
(v) Let $\A$ be a filter, then $\bigcap\overline{\A}$ is a pure
filter and $\bigcap\overline{\A}\subseteq\A$. In addition, if $\A$
is a pure filter then $\A=\bigcap\overline{\A}$.\\
(vi) Let $\A$ and $\mathcal{B}$ be pure filters then $\A\subseteq
\mathcal{B}$ if and only if
$\overline{\mathcal{B}}\subseteq\overline{\A}$, also $\A=\B$ if and
only if $\overline{\A}=\overline{\B}$.
\end{lem}
\textbf{Proof :} (i), (ii) and (iii) are clear.

 (iv) It is obvious
that $\A=\bigcap\mathcal{J}$ is a pure filter and
$\mathcal{J}\subseteq \overline{\A}$. Let $f\in\f$ such that
$\mathcal{J}\subseteq (\overline{Z(f)})^{\circ}$. Since for each
$\widetilde{p}\in\mathcal{J}\subseteq (\overline{Z(f)})^{\circ}$ we
have $Z(f)\in\widetilde{p}$ (Lemma 2.2), so $Z(f)\in\A$. Hence by
$(i)$
\[
\overline{\A}=\bigcap_{A\in\A}\overline{A}\subseteq\bigcap_{\mathcal{J}\subseteq
(\overline{A})^{\circ}}\overline{A}\subseteq\mathcal{J},
\]
and therefore $\overline{\A}=\mathcal{J}$.

 (v) It is obvious that
$\bigcap\overline{\A}$ is a pure filter and
$\bigcap\overline{\A}\subseteq\A$. Now let $\A$ be a pure filter,
then for every $\widetilde{p}\in\overline{\A}$,
$\A\subseteq\widetilde{p}$ and this implies that
$\A\subseteq\bigcap_{\widetilde{p}\in\overline{\A}}\widetilde{p}
=\bigcap\overline{\A}$.

 (vi) It is clear.$\hfill\Box$
\begin{defn}
Let $A$ be a nonempty subset of $\mathcal{R}$.  We call $\bigcap A$
the pure filter generated by $A$.
\end{defn}
\begin{defn}
Let $\A,\B\subseteq Z(\f)$ be filters. We define
\begin{eqnarray*}
  \A \odot\B &=& \{A\in
Z(\f): \,\,Z(\mathbf{T}_{\mu}f)\in\A ,\,\,\forall f\in
[A],\,\,\forall \A^{\mu}\in\overline{\B}\},
 \end{eqnarray*}
where $[A]=\{f\in \f^{+}: Z(f)=A\}$
\end{defn}
\begin{lem}
Let $\A, \B$ be filters then $\A \odot\B$ is a filter.
\end{lem}
\textbf{Proof :} Let $A, B\in\A \odot\B$ and
$\A^{\mu}\in\overline{\B}$, then for each $h\in [A\cap B]$ there
exist $u\in [A]$ and $v\in [B]$ such that $0\leq h\leq u+v$ and so
for each $\mu\in S^{\f}$,
$\mathbf{T}_{\mu}h\leq\mathbf{T}_{\mu}(u+v)$. Therefore
$Z(\mathbf{T}_{\mu}u)\cap
Z(\mathbf{T}_{\mu}v)=Z(\mathbf{T}_{\mu}(u+v))\subseteq
Z(\mathbf{T}_{\mu}h)$ implies that $Z(\mathbf{T}_{\mu}h)\in\A$ for
each $h\in [A\cap B]$ and $\A^{\mu}\in \overline{\mathcal{B}}$.
 Hence  $A\cap B\in\A \odot\B$.\\
Now let $A\in \A \odot\B$ and $B\in Z(\f)$ such that $A\subseteq B$.
Choose $\A^{\mu}\in\overline{\B}$. For each $h\in [B]$, there exists
$f\in [A]$ such that $0\leq h\leq f$. Therefore
$Z(\mathbf{T}_{\mu}f)\subseteq Z(\mathbf{T}_{\mu}h)$ and this
implies $B\in \A\odot\B$.$\hfill\Box$
\begin{defn}
Let $\A$ and $\B$ be filters, we define  $\A+ \B=\bigcap\overline{\A
\odot\mathcal{B}}$.
\end{defn}
By Definition 3.4, $\A+ \B$ is a pure filter generated by
$\overline{\A\odot\B}$. The $"+"$ operation on $\mathcal{P}(\f)$
 is not necessarily commutative.
\begin{lem}
Let $\mathcal{A}$ and $\mathcal{B}$ be pure filters and $\mu , \nu
\in S^{\mathcal{F}}$, then\\
(i) $\mathcal{A}+\mathcal{A}^{\mu}=\bigcap_{\mathcal{A}^{\nu}\in
\overline{\mathcal{A}}}\mathcal{A}^{\nu\mu},$\\
(ii)
$\mathcal{A}^{\mu}+\mathcal{A}^{\nu}=\mathcal{A}^{\mu\nu}=\mathcal{A}^{\mu}*\mathcal{A}^{\nu},$\\
(iii) $\mathcal{A}+
\mathcal{B}\subseteq\bigcap_{\mathcal{A}^{\nu}\in\overline{\mathcal{B}},
\mathcal{A}^{\mu}\in\overline{\mathcal{A}}} \mathcal{A}^{\mu\nu},$ \\
(iv) $\mathcal{A}^{\mu}+
\mathcal{A}\subseteq\bigcap_{\mathcal{A}^{\nu}\in
\overline{\mathcal{A}}}\mathcal{A}^{\mu\nu}.$
\end{lem}
\textbf{Proof :} $(i)$ Let $D=\{\nu\mu:\nu\in\overline{\A}\}$ and
$\eta\notin D$, then there exists $f\in\f$ such that $f:S\rightarrow
[0,1]$, $D\subseteq (\overline{Z(f)})^{\circ}$ and
$\eta\notin\overline{Z(f)}$. Therefore
$\nu\in\rho^{-1}_{\mu}(D)\subseteq\rho^{-1}_{\mu}((\overline{Z(f)})^{\circ})\subseteq
\overline{Z(\mathbf{T}_{\mu}g)}$ for every $g\in [Z(f)]$ and so
$Z(\mathbf{T}_{\mu}g)\in\A$ for every $f\in [Z(f)]$. Hence
$Z(f)\in\A\odot\A^{\mu}$, and this implies that
$\eta\notin\overline{\A\odot\A^{\mu}}$. So
$\overline{\A\odot\A^{\mu}}\subseteq D$ and we will have
$\bigcap_{\mathcal{A}^{\nu}\in
\overline{\mathcal{A}}}\mathcal{A}^{\nu\mu}\subseteq\mathcal{A}+\mathcal{A}^{\mu}.$\\
Now let $D=\{\nu\mu:\nu\in\overline{\A}\}$, if $A\in\A\odot\A^{\mu}$
then $Z(\mathbf{T}_{\mu}f)\in\A^{\nu}$ for every $f\in [A]$ and
$\A^{\nu}\in\overline{\A}$. Hence
$\widehat{f}(\nu\mu)=\nu\mu(f)=\nu(\mathbf{T}_{\mu}f)=0$ and so
$\nu\mu\in\overline{A}$. Therefore we will have
$D\subseteq\bigcap_{A\in\A\odot\A^{\mu}}\overline{A}=\overline{\A\odot\A^{\mu}}$,
and this completes the proof.

$(ii)$ Obvious.

 $(iii)$ Let $A\in\A\odot\B$ then
$Z(\mathbf{T}_{\mu}f)\in\A^{\mu}$ for each
$\A^{\nu}\in\overline{\B}$, $\A^{\mu}\in\overline{\A}$ and $f\in
[A]$. Therefore
$\widehat{f}(\nu\mu)=\nu\mu(f)=\nu(\mathbf{T}_{\mu}f)=0$ for each
$\A^{\nu}\in\overline{\B}$, $\A^{\mu}\in\overline{\A}$ and $f\in
[A]$ so $\nu\mu\in\overline{A}$, then we will have
$$
\{\mu\nu\in
S^{\f}:\A^{\mu}\in\overline{\A},\,\,\A^{\nu}\in\overline{\B}\}
\subseteq\bigcap_{A\in\A\odot\B}\overline{A}=\overline{\A\odot\B},
$$
and this implies that $\mathcal{A}+
\mathcal{B}=\bigcap\overline{\A\odot\B}\subseteq
\bigcap_{\mathcal{A}^{\nu}\in\overline{\mathcal{B}},
\mathcal{A}^{\mu}\in\overline{\mathcal{A}}} \mathcal{A}^{\mu\nu}.$

 $(iv)$ Obvious.$\hfill\Box$
\begin{lem}
Let $\mathcal{A}$ be a pure filter such that
$\mathcal{A}\subseteq\mathcal{A}+\mathcal{A}$. Then
$\overline{\mathcal{A}}$ is a subsemigroup of $\mathcal{R}$.
\end{lem}
\textbf{Proof :} It is obvious that
$\A\subseteq\mathcal{A}+\mathcal{A}\subseteq
\bigcap_{\mathcal{A}^{\mu},\mathcal{A}^{\nu}
\in\overline{\mathcal{A}}}\mathcal{A}^{\mu\nu}$. So
 $\mathcal{A}^{\mu},\mathcal{A}^{\nu} \in\overline{\mathcal{A}}$
 implies that $\mathcal{A}\subseteq\mathcal{A}^{\mu\nu}$ and so
$\mathcal{A}^{\mu\nu}\in\overline{\mathcal{A}}$. Hence
$\overline{\A}$ is a subsemigroup.$\hfill\Box$
\begin{rem}
If $\A$ is a filter such that $\A\subseteq\A\odot\A$, thenthe conclusion of Lemma 3.9
remains true since
$\overline{\A+\A}=\overline{\A\odot\A}\subseteq\overline{\A}$.
\end{rem}
\begin{defn}
Let $\A$ be a filter and $A\in Z(\f)$; we define
\begin{eqnarray*}
\Omega_{\A}(A)&=&\{x\in
S:\overline{\A}\subseteq\overline{\lambda_{x}^{-1}(A)}\},\\
\Omega_{\A^{\circ}}(A)&=&\{x\in S: \overline{\A}\subseteq
 \lambda_{\widehat{x}}^{-1}((\overline{A})^{\circ})\}.
 \end{eqnarray*}
\end{defn}
Clearly, for each $A\in Z(\f)$ and filter $\A$ we have
$\Omega_{\mathcal{A}^{\circ}}(A)\subseteq\Omega_{\mathcal{A}}(A)$.
\begin{lem}
Let $\A$ and $\B$ be filters and $A, \,\, B\in Z(\f)$. Then

 (i)
$\Omega_{\A}(A)=\bigcap_{\mu\in\overline{\mathcal{A}},\,f\in
[A]}Z(\mathbf{T}_{\mu}f)$,

(ii) $\Omega_{\A}(A\cap B)=\Omega_{\A}(A)\cap\Omega_{\A}(B)$,

 (iii)
If $\A\subseteq\B$ then $\Omega_{\A}(A)\subseteq\Omega_{\B}(A)$,

(iv) $A\subseteq B$ then $\Omega_{\A}(A)\subseteq\Omega_{\A}(B)$,

(v) For each $x\in S$,
$\lambda_{x}^{-1}(\Omega_{\A}(A))=\Omega_{\A}(\lambda_{x}^{-1}(A))$,

(vi) For all $x\in \Omega_{\A^{\circ}}(A)$ then
$\lambda_{x}^{-1}(A)\in\A$.
\end{lem}
\textbf{Proof :} $(i)$ Clearly we have:
\begin{eqnarray*}
\Omega_{\mathcal{A}}(A)&=&\{x\in
S:\overline{\A}\subseteq\overline{\lambda_{x}^{-1}(A)}\}\\
&=&\{x\in S:\overline{\A}\subseteq
\lambda_{\widehat{x}}^{-1}(\overline{A})\}\\
&=&\{x\in
S:\widehat{x}*\overline{\A}\subseteq\overline{Z(f)},\,\,\,\,\forall
f\in
[A]\}\\
&=&\{x\in S:
\varepsilon(x)\mu(f)=0\,\,\,\mu\in\overline{\A},\,\,\,\,\forall f\in [A]\}\\
&=& \{x\in S:
\mathbf{T}_{\mu}(f)(x)=0,\,\,\,\forall\mu\in\overline{\A},\,\,\,\,\forall f\in [A]\}\\
&=&\bigcap_{\mu\in\overline{\mathcal{A}},\,f\in
[A]}Z(\mathbf{T}_{\mu}f).
\end{eqnarray*}

$(ii)$ It is obvious that  $[A]+[B]\subseteq [A\cap B]$ for all $A,
\,\, B\in Z(\f)$. Also for each $f\in [A\cap B]$ there exist
$g_{_{f}}\in [A]$ and $h_{f}\in [B]$ such that $f\leq
g_{_{f}}+h_{_{f}}$, and so $\mathbf{T}_{\mu}f\leq
\mathbf{T}_{\mu}(g_{_{f}}+h_{_{f}})$ for all $\mu\in S^{\f}$, so
$Z(\mathbf{T}_{\mu}(g_{_{f}}+h_{_{f}}))\subseteq
Z(\mathbf{T}_{\mu}(f))$. Now for any filter $\A$, we can write
\begin{eqnarray*}
\Omega_{\A}(A)\cap\Omega_{\A}(B)
&=&(\bigcap_{\mu\in\overline{\mathcal{A}},\,g_{_{f}}\in
[A]}Z(\mathbf{T}_{\mu}g_{f}))\cap(\bigcap_{\mu\in\overline{\mathcal{A}},\,h_{f}\in
[B]}Z(\mathbf{T}_{\mu}h_{_{f}}))\\
&\subseteq & \bigcap_{\mu\in\overline{\mathcal{A}},\,g_{_{f}}\in
[A],\,\,h_{_{f}}\in
[B]}(Z(\mathbf{T}_{\mu}g_{_{f}})\cap Z(\mathbf{T}_{\mu}h_{_{f}}))\\
&=&\bigcap_{\mu\in\overline{\mathcal{A}},\,g_{_{f}}\in
[A],\,\,h_{_{f}}\in
[B]}(Z(\mathbf{T}_{\mu}(g_{_{f}}+h_{_{f}}))\\
&\subseteq &\bigcap_{\mu\in\overline{\mathcal{A}},\,f\in
[A\cap B]}Z(\mathbf{T}_{\mu}f)\\
&=& \Omega_{\A}(A\cap B)\\
&=&\bigcap_{\mu\in\overline{\mathcal{A}},\,f\in
[A\cap B]}Z(\mathbf{T}_{\mu}f)\\
&\subseteq&\bigcap_{\mu\in\overline{\A}}(\bigcap_{g_{_{f}}\in
[A],\,\,h_{_{f}}\in
[B]}Z(\mathbf{T}_{\mu}(g_{_{f}}+h_{_{f}})))\\
&=&\Omega_{\A}(A)\cap\Omega_{\A}(B).
\end{eqnarray*}

$(iii)$ Since $\A\subseteq\B$ therefore
$\overline{\B}\subseteq\overline{\A}$, and so
\begin{eqnarray*}
\hspace{3cm}\Omega_{\A}(A)&=& \{x\in S: \overline{\A}\subseteq
 \overline{\lambda_{x}^{-1}(A)}\}\\
 &=& \{x\in S: \overline{\B}\subseteq\overline{\A}\subseteq
 \overline{\lambda_{x}^{-1}(A)}\}\\
 &\subseteq & \{x\in S: \overline{\B}\subseteq
 \overline{\lambda_{x}^{-1}(A)}\}\\
 &=& \Omega_{\B}(A).
\end{eqnarray*}

$(iv)$ Obvious.

$(v)$ For every $x\in S$ we will have
\begin{eqnarray*}
\lambda_{x}^{-1}(\Omega_{\mathcal{B}}(B))&=&\lambda_{x}^{-1}(\{t\in
S:\overline{\mathcal{B}}\subseteq\overline{\lambda_{t}^{-1}(B)}\})\\
&=&\{y\in
S:\overline{\mathcal{B}}\subseteq\overline{\lambda_{xy}^{-1}(B)}\}\\
 &=&\{y\in
S:\overline{\mathcal{B}}\subseteq\overline{\lambda_{y}^{-1}(\lambda_{x}^{-1}(B))}\}\\
&=&\Omega_{\mathcal{B}}(\lambda_{x}^{-1}(B)).
\end{eqnarray*}

$(vi)$ For every $x\in\Omega_{\mathcal{A}^{\circ}}(A)$ we will have
\begin{eqnarray*}
 \overline{\A} &\subseteq & \lambda_{\widehat{x}}^{-1}((\overline{A})^{\circ})\\
&\subseteq & (\lambda_{\widehat{x}}^{-1}(\overline{A}))^{\circ}\\
&\subseteq & \overline{\lambda_{x}^{-1}(A)},
\end{eqnarray*}
and this implies that $\lambda_{x}^{-1}(A)\in\A $, (Lemma
2.2.).$\hfill\Box$

\begin{thm}
 Let $\mathcal{A}$ be a filter.
 Then the following statements are equivalent.\\
(i) $\overline{\mathcal{A}}$ is a left ideal of $\mathcal{R}$,\\
(ii) For all $A\in\A,\hspace{.2cm}\Omega_{\A}(A)=S$,\\
 (iii) For all $A\in Z(\f)$, if $\,\,\overline{\mathcal{A}}\subseteq
(\overline{A})^{\circ}$ then  $\Omega_{\A^{\circ}}(A)=S$.
\end{thm}
\textbf{Proof :} $(ii)\rightarrow (i)$ Let  $\Omega_{\A}(A)=S$ for
all $A\in\A$, so $Z(\mathbf{T}_{\mu}f)=S$ for each
$\A^{\mu}\in\overline{\A}$ and $f\in [A]$. Hence
$Z(\mathbf{T}_{\mu}f)=S\in\A^{\nu}$ for each
$\A^{\nu}\in\mathcal{R}$, $\A^{\mu}\in\overline{\A}$ and $f\in [A]$.
 Therefore $\widehat{f}(\nu\mu)=\nu\mu(f)=\nu(\mathbf{T}_{\mu}f)=0$ for
 all $\A^{\nu}\in\mathcal{R}$, $\A^{\mu}\in\overline{\A}$ and $f\in [A]$ and so
 $\nu\mu\in \overline{A}$. This implies that
$\mathcal{R}*\overline{\A}\subseteq\bigcap_{A\in\A}\overline{A}=\overline{\A}$.

$(i)\rightarrow (ii)$ Let $\overline{\A}$ be a left ideal, then
$\widehat{x}*\overline{\A}\subseteq\overline{\A}\subseteq\overline{A}$
for each $A\in\A$ and $x\in S$. So
$\overline{\A}\subseteq\overline{\lambda_{x}^{-1}(A)}$ for all $A\in
\A$ and $x\in S$ and then implies that $\forall A\in\A,
\hspace{.2cm}\Omega_{\A}(A)=S$.

 $(iii)\rightarrow (i)$
 Let for each $A\in Z(\f)$,
$\,\,\,\overline{\mathcal{A}}\subseteq(\overline{A})^{\circ}$
implies that $\Omega_{\mathcal{A}^{\circ}}(A)=S$, since
$\Omega_{\A^{\circ}}(A)\subseteq \Omega_{\A}(A)=\bigcap_{f\in [A],\,
\mu\in \overline{\A}}Z(\mathbf{T}_{\mu}f)$, by Lemma 3.12, we will
have $A\in \A^{\nu}+\A^{\mu}=\A^{\nu\mu}$ for each $\A^{\nu}\in
\mathcal{R}$ and $\A^{\mu}\in \overline{\A}$. Now
$\overline{\A}=\bigcap_{A\in
\Gamma}(\overline{A})^{\circ}\subseteq\bigcap_{A\in
\Gamma}\overline{A}$, where $\Gamma=\{A\in
Z(\f):\overline{\A}\subseteq(\overline{A})^{\circ}\}$, implies that
$\mathcal{R}* \overline{\A}\subseteq \overline{\A}$. Hence
$\overline{\A}$
 is a left ideal.

$(i)\rightarrow (iii)$ Let $A\in Z(\f)$ such that
$\overline{\A}\subseteq(\overline{A})^{\circ}$, since
$\overline{\A}$ is a left ideal then $\mathcal{R}*
\overline{\A}\subseteq \overline{\A}\subseteq
(\overline{A})^{\circ}$. Therefore for all $\A^{\mu}\in \mathcal{R}$
and $\A^{\nu}\in \overline{\A}$ we have $\A^{\mu\nu}\in
(\overline{A})^{\circ}$, and so for each $x\in S$,
 $\varepsilon(x)\nu\in (\overline{A})^{\circ}$. Hence $\nu\in
\lambda_{\widehat{x}}^{-1}((\overline{A})^{\circ})$ for all
$\A^{\nu}\in \overline{\A}$ and $x\in S$, so $
\Omega_{\A^{\circ}}(A)=\{x\in S: \overline{\A}\subseteq
\lambda_{\widehat{x}}^{-1}((\overline{A})^{\circ})\}=S$.$\hfill\Box$
\begin{cor}
Let $\A$ be a pure filter, then $\overline{\A}$ is a left ideal if
and only if $\forall A\in\A,\hspace{.2cm}\Omega_{\A}(A)=S$.
\end{cor}
\textbf{Proof :} It is obvious.$\hfill\Box$

 Therefore if $\A$ is a
pure filter, we will have $"\forall
A\in\A,\,\,\,\,\Omega_{\A}(A)=S"$ if and only if $\overline{\A}$ is
a left ideal. Since for each $\A, \mathcal{B}\in \mathcal{P}(\f)$,
$\A\subseteq\mathcal{B}$ if and only if
$\overline{\mathcal{B}}\subseteq\overline{\A}$, then for $\A\in
\mathcal{P}(\f)$, $\overline{\A}$ is a minimal left ideal in
$\mathcal{R}$ if and only if $\A$ is maximal in $\mathcal{P}(\f)$
with respect to the property $"\forall
A\in\A,\,\,\,\Omega_{\A}(A)=S"$. This concept is similar to the
discrete case, (see \cite{Talin1}).
\begin{thm}
Let $\mathcal{A}$ be a pure filter, then $\overline{\mathcal{A}}$ is
a right ideal of $\mathcal{R}$ if and only if
$\mathcal{A}\subseteq\mathcal{A}+\mathcal{A}^{\nu}$ for any
$\mathcal{A}^{\nu}\in\mathcal{R}$.
\end{thm}
\textbf{Proof :} If $\overline{\mathcal{A}}$ is a right ideal of
$\mathcal{R}$, then
$\overline{\mathcal{A}}*\mathcal{R}\subseteq\overline{\mathcal{A}}$,
and so $\mathcal{A}\subseteq\mathcal{A}^{\mu}+\mathcal{A}^{\nu}$ for
any $\mathcal{A}^{\mu}\in\overline{\mathcal{A}}$ and
$\mathcal{A}^{\nu}\in\mathcal{R}$. Then for $A\in\mathcal{A}$ we
have $Z(\mathbf{T}_{\nu}f)\in\mathcal{A}^{\mu}$ for all $f\in [A]$,
$\mathcal{A}^{\nu}\in\mathcal{R}$ and
$\mathcal{A}^{\mu}\in\overline{\mathcal{A}}$. Therefore
$Z(\mathbf{T}_{\nu}f)\in\mathcal{A}$ for all $f\in [A]$ and
$\mathcal{A}^{\nu}\in\mathcal{R}$, and so
$A\in\mathcal{A}+\mathcal{A}^{\nu}$ for all
$\mathcal{A}^{\nu}\in\mathcal{R}$. Hence
$\mathcal{A}\subseteq\mathcal{A}+\mathcal{A}^{\nu}$ for all
$\mathcal{A}^{\nu}\in\mathcal{R}$.

 Conversely, if
$\mathcal{A}\subseteq\mathcal{A}+\mathcal{A}^{\nu}$ for all
$\mathcal{A}^{\nu}\in\mathcal{R}$, then
$\mathcal{A}\subseteq\mathcal{A}^{\mu}+\mathcal{A}^{\nu}
=\mathcal{A}^{\mu\nu}$ for each $\mathcal{A}^{\nu}\in\mathcal{R}$
and $\mathcal{A}^{\mu}\in\overline{\mathcal{A}}$. Hence
$\overline{\mathcal{A}}*\mathcal{R}\subseteq\overline{\mathcal{A}}$.$\hfill\Box$
\begin{cor}
Let $S$ be a commutative semigroup. Then for every pure filter
$\mathcal{A}$, $\overline{\mathcal{A}}$ is an ideal of $\mathcal{R}$
if and only if $\A\subseteq\A+\A^{\nu}$ for every
$\A^{\nu}\in\mathcal{R}$.
\end{cor}
\textbf{Proof :} If $\A\subseteq\A+\A^{\nu}$ for every
$\A^{\nu}\in\mathcal{R}$, then $\overline{\A}$ is a right ideal by
Theorem 3.15.

Conversely, since $S$ is a commutative semigroup and
$\mathcal{A}\subseteq\mathcal{A}+\mathcal{A}^{\varepsilon(x)}$ for
any $x\in S$, therefore $A\in\mathcal{A}$ implies
$Z(L_{x}f)=Z(\mathbf{T}_{\varepsilon(x)}f)\in\mathcal{A}$ for each
$f\in [A]$ and  $x\in S$. Hence we will have
$Z(\mathbf{T}_{\mu}f)=S$ for all $f\in [A]$ and
$\mathcal{A}^{\mu}\in \overline{\mathcal{A}}$, and so Theorem 3.13
implies $\overline{\mathcal{A}}$ is a left ideal.$\hfill\Box$
\begin{lem}
Let $\A$ be a filter such that $\overline{\A}$ is a left ideal of
$\mathcal{R}$. Then $\mathcal{D}_{\mathcal{B}}=\{B\in
Z(\f):\Omega_{\mathcal{B}}(B)=S\}$ is a
 filter for any filter $\A\subseteq\mathcal{B}$,
  and $\overline{\mathcal{D}_{\mathcal{B}}}$ is a left ideal.
  Further, if $\overline{\A}$ is a minimal left ideal then
  $\overline{\mathcal{D}_{\B}}$ is a minimal left ideal.
\end{lem}
\textbf{Proof :} Let $\mathcal{B}$ be a filter such that
$\A\subseteq\mathcal{B}$ and $\mathcal{D}_{\B}=\{B\in
Z(\f):\Omega_{\mathcal{B}}(B)=S\}$. By Lemma 3.12,
$\mathcal{D}_{\B}$ is a filter because
$$\Omega_{\mathcal{B}}(A\cap
B)=\Omega_{\mathcal{B}}(A)\cap\Omega_{\mathcal{B}}(B)=S.$$ Also for
all $A,B\in Z(\f)$, $A\subseteq B$ implies that
$\Omega_{\mathcal{B}}(A)\subseteq\Omega_{\mathcal{B}}(B)$. Now let
$\overline{\A}$ is a left ideal of $\mathcal{R}$ therefore $\A$ has
the property $"\forall A\in\A,\;\;\Omega_{\A}(A)=S"$, (Theorem
3.13). By Lemma 3.12,
$S=\Omega_{\A}(A)\subseteq\Omega_{\mathcal{B}}(A)$ for all $A\in\A$
 and this implies $\A\subseteq\mathcal{D}_{\B}$. Take
 $B\in\mathcal{D}_{\B}$, then $\Omega_{\mathcal{B}}(B)=S$, and so
$S=\lambda_{x}^{-1}(\Omega_{\mathcal{B}}(B))=
\Omega_{\mathcal{B}}(\lambda_{x}^{-1}(B))$ for all $x\in S$, (Lemma
3.12). Hence $\lambda_{x}^{-1}(B)\in \mathcal{D}_{\B}$ for each
$x\in S$, and this implies $\Omega_{\mathcal{D}_{\B}}(B)=S$.
Therefore $\overline{\mathcal{D}_{\B}}$ is a left ideal, by Theorem
3.13.

 Now let $\overline{\A}$ be a minimal left ideal of
$\mathcal{R}$. Since $\A\subseteq \mathcal{D}_{\B}$ so
$\overline{\mathcal{D}_{\B}}\subseteq\overline{\A}$. But
$\overline{\A}$ is a minimal left ideal and
$\overline{\mathcal{D}_{\B}}$ is a left ideal, therefore
$\overline{\mathcal{D}_{\B}}=\overline{\A}$ and this implies
  $\overline{\mathcal{D}_{\B}}$ is a minimal left ideal.$\hfill\Box$
\begin{thm}
Let $\A$ be a pure filter such that $\overline{\A}$ is a left ideal
of $\mathcal{R}$. Then the following statements are equivalent.\\
(i) $\overline{\A}$ is a minimal left ideal.\\
 (ii) $\mathcal{D}_{\mathcal{B}}=\{B\in Z(\f):\Omega_{\mathcal{B}}(B)=S\}$
is a filter for every filter $\A\subseteq\mathcal{B}$ and
 $\overline{\mathcal{D}_{\mathcal{B}}}$ generates $\A$ as a pure
 filter.\\
(iii) $\LL_{\B}=\{B\in\mathcal{B}:\Omega_{\mathcal{B}}(B)=S\}$ is a
filter for every filter $\A\subseteq\mathcal{B}$ and
$\overline{\LL_{\B}}$ generates $\A$ as a pure filter.
\end{thm}
\textbf{Proof :} $(i)\rightarrow (ii)$

 By Lemma 3.17, it is obvious
that $\mathcal{D}_{\B}$ is a filter and also
$\overline{\mathcal{D}_{\B}}=\overline{\A}$ is a minimal left ideal.
Hence $\mathcal{D}_{\B}$ generates $\A$ as a pure filter.

$(ii)\rightarrow (iii)$ For every  $\A\subseteq\mathcal{B}$,
$\LL_{\B}=\{B\in\B:\Omega_{\mathcal{B}}(B)=S\}$ is a filter, similar
to Lemma 3.17. Now since $\A\subseteq\LL_{\B}=\{B\in
\mathcal{B}:\Omega_{\mathcal{B}}(B)=S\}\subseteq\{B\in
Z(\f):\Omega_{\mathcal{B}}(B)=S\}$ so the result is immediate.

$(iii)\rightarrow (i)$ Let $L$ be a minimal left ideal of
$\mathcal{R}$ such that $L\subseteq\overline{\A}$. Let $\LL=\cap L$,
so $\LL$ is a pure filter and $\A\subseteq\LL$. Therefore
$\Omega_{\mathcal{L}}(B)=S$ for all $B\in\mathcal{L}$, and so
$\overline{\A}\subseteq \overline{B}$. Hence
$\overline{\A}\subseteq\bigcap_{B\in\LL}\overline{B}=\overline{\LL}$
and this implies $\overline{\A}$ is a minimal left ideal of
$\mathcal{R}$.$\hfill\Box$
\begin{cor}
Let $\A$ be a filter such that $\overline{\A}$ is a left ideal of
$\mathcal{R}$. If $\overline{\A}$ is a minimal left ideal of
$\mathcal{R}$, then for any $A\in Z(\f)$, $\Omega_{\A}(A)=S$ implies
$\overline{\A}\subseteq\overline{A}$.
\end{cor}
\textbf{Proof :} Take $\mathcal{B}=\A$ in Theorem 3.18 $(ii)$. The
result is immediate.$\hfill\Box$
\begin{thm}
Let $\A$ be a pure filter such that $\overline{\A}$ is a left ideal
of $\mathcal{R}$. Then the following statements are equivalent.\\
(i) $\overline{\A}$ is a minimal left ideal of $\mathcal{R}$.\\
 (ii) For $A\in Z(\f)$, if $\overline{\A}-\overline{A}\neq\emptyset$ then
for each  $B\in Z(\f)$ with the property
$\overline{S}-\overline{A}\subseteq
  (\overline{B})^{\circ}$, there exists $F\in [S]^{<\omega}$ such that $\bigcup_{x\in
F}\lambda_{x}^{-1}(B)\in\A$.\\
 (iii) Let $A\in Z(\f)$, if
$\overline{\A}-\overline{A}\neq\emptyset$ then for each $B\in Z(\f)$
with the property $\overline{S}-\overline{A}\subseteq
(\overline{B})^{\circ}$, there exist $C\in\A$ and $F\in
[S]^{<\omega}$ such that $C\cap A\subseteq\bigcup_{x\in
F}\lambda_{x}^{-1}(B)$.
\end{thm}
\textbf{Proof :} $(i)\rightarrow (ii)$ Assume that $\overline{\A}$
is a minimal left ideal of $\mathcal{R}$ and suppose that ${(ii)}$
is not true. Assume further that there exists $A\in Z(\f)$ such that
$\overline{\A}-\overline{A}\neq\emptyset$ and for some $B\in Z(\f)$
with the property  $\overline{S}-\overline{A}\subseteq
(\overline{B})^{\circ}$ and for all $F\in [S]^{<\omega}$ implies
$\bigcup_{x\in F}\lambda_{x}^{-1}(B)\notin\A$. Now let
$\Lambda=\{D\cap(\bigcap_{x\in F}\lambda_{x}^{-1}(A)):
\overline{\A}\subseteq (\overline{D})^{\circ},\,\,D\in
Z(\f),\;\;F\in [S]^{<\omega}\}$. If there exists $x\in S$ such that
$\lambda_{x}^{-1}(A)=\emptyset$ since $A\cup B=S$, we will have
$\lambda_{x}^{-1}(B)\in\A$ and this is a contradiction. Therefore
$\lambda_{x}^{-1}(A)\neq\emptyset$ for every $x\in S$. Now if there
exists $F\in [S]^{<\omega}$ and $D\in\A$ such that
$D\cap(\bigcap_{x\in F}\lambda_{x}^{-1}(A))=\emptyset$, then $A\cup
B=S$ implies $ D \subseteq\bigcup_{x\in F}\lambda_{x}^{-1}(B)$ and
this is a contradiction. So $\Lambda$ has the finite intersection
property and generates a filter $\mathcal{L}$. Since $\overline{\A}$
is a left ideal therefore for every $D\in Z(\f)$ such that
$\overline{\A}\subseteq(\overline{D})^{\circ}$ then we will have
$\overline{\A}\subseteq\lambda_{\widehat{x}}^{-1}((\overline{D})^{\circ})
\subseteq(\overline{\lambda_{x}^{-1}(D)})^{\circ}$ for each $x\in
S$, and so  $\lambda_{x}^{-1}(D)\in\LL$   for every $x\in S$.
Therefore $\Omega_{\mathcal{L}}(E)=S$  for every $E\in\mathcal{L}$
and this implies $\overline{\LL}$ is left ideal, by Theorem 3.10.

Since $
 \overline{\LL}=\bigcap_{E\in\LL}\overline{E}\subseteq
\bigcap_{D\in\A,\,\,\overline{\A}\subseteq(\overline{D})^{\circ}}\overline{D}
=\overline{\A},$ and $\overline{\A}$ is a minimal left ideal,
therefore $\overline{\A}=\overline{\LL}$ and $\overline{\LL}$ is a
minimal left ideal. Also
$\LL\subseteq\mathcal{D}_{\mathcal{L}}=\{B\in
Z(\f):\Omega_{\mathcal{L}}(B)=S\}$ and so $\overline{D}_{\LL}$ is a
minimal left ideal, by Lemma 3.15. It is obvious that $A\in
\mathcal{D}_{\LL}$ and this implies
$\overline{\A}=\overline{\mathcal{D}_{\LL}}\subseteq\overline{A}$.
So we have a contradiction.

$(ii)\rightarrow (iii)$ Put $C=\bigcup_{x\in
F}\lambda_{x}^{-1}(B)\in\A$, then $C\cap A\subseteq\bigcup_{x\in
F}\lambda_{x}^{-1}(B)$.

$(iii)\rightarrow (i)$ Let for each $A\in Z(\f)$ with the property
$\overline{\A}-\overline{A}\neq\emptyset$ and for each $B\in Z(\f)$
that $\overline{S}-\overline{A}\subseteq (\overline{B})^{\circ}$
there exist $C\in\A$ and $F\in [S]^{<\omega}$ such that $C\cap
A\subseteq\bigcup_{x\in F}\lambda_{x}^{-1}(B)$. Suppose that
$\overline{\A}$ is not a minimal left ideal, therefore there exists
a pure filter $\mathcal{U}$ such that $\A\subseteq\mathcal{U}$ and
$\overline{\mathcal{U}}$ be a minimal left ideal. So there exists
$B_{\circ}\in \mathcal{U}$ such that
$\overline{\mathcal{U}}\subseteq (\overline{B_{\circ}})^{\circ}$ and
$\overline{\A}-\overline{B_{\circ}}\neq\emptyset$. Hence for each
$B\in Z(\f)$ such that $\overline{S}-\overline{B_{\circ}}\subseteq
(\overline{B})^{\circ}$, there exist $C\in\A$ and $F_{B}\in
[S]^{<\omega}$ such that $C\cap B_{\circ}\subseteq\bigcup_{x\in
F_{B}}\lambda_{x}^{-1}(B)$. Therefore $\A\subseteq\mathcal{U}$
implies $C\cap B_{\circ}\in\mathcal{U}$, and this conclude
$\bigcup_{x\in F_{B}}\lambda_{x}^{-1}(B)\in\mathcal{U}$. Since
$\bigcap_{x\in F_{B}}\lambda_{x}^{-1}(B_{\circ})\in\mathcal{U}$ so
for each $B\in Z(\f)$ that
$\overline{S}-\overline{B_{\circ}}\subseteq (\overline{B})^{\circ}$
then there exists $F_{B}\in [S]^{<\omega}$ such that
\begin{eqnarray*}
(\bigcup_{x\in F_{B}}\lambda_{x}^{-1}(B))\cap (\bigcap_{x\in
F_{B}}\lambda_{x}^{-1}(B_{\circ}))&=&\bigcup_{x\in F_{B}}
(\lambda_{x}^{-1}(B)\cap\bigcap_{y\in
F_{B}}\lambda_{y}^{-1}(B_{\circ}))\\
&\subseteq&\bigcup_{y\in F_{B}}\lambda_{y}^{-1}(B\cap B_{\circ})\\
&\in &\mathcal{U},
\end{eqnarray*}
and so
\begin{eqnarray*}
\overline{\mathcal{U}}&\subseteq&\bigcup_{y\in
F_{B}}\lambda_{\widehat{y}}^{-1}(\overline{B\cap B_{\circ}})\\
&\subseteq&\bigcup_{y\in
F_{B}}\lambda_{\widehat{y}}^{-1}(\overline{B}\cap
\overline{B_{\circ}})\\
&\subseteq&\bigcup_{y\in\,\,
S}\lambda_{\widehat{y}}^{-1}(\overline{B}\cap \overline{B_{\circ}}).
\end{eqnarray*}
Since $\overline{\mathcal{U}}\subseteq\bigcup_{y\in
F_{B}}\lambda_{\widehat{y}}^{-1}((\overline{B_{\circ}})^{\circ})
\subseteq\bigcup_{y\in
S}\lambda_{\widehat{y}}^{-1}((\overline{B_{\circ}})^{\circ})$ so we
have $\overline{\mathcal{U}}\subseteq\bigcup_{y\in
S}\lambda_{\widehat{y}}^{-1}(\overline{B}\cap
(\overline{B_{\circ}})^{\circ})$ for every $B\in Z(\f)$ that
$\overline{S}-\overline{B_{\circ}}\subseteq (\overline{B})^{\circ}$,
and so
$$\overline{\mathcal{U}}\subseteq\bigcup_{y\in
S}\lambda_{\widehat{y}}^{-1}(\overline{\overline{B_{\circ}}^{c}}\cap
(\overline{B_{\circ}})^{\circ})=\emptyset.$$
 Now we have a contradiction.$\hfill\Box$
\bibliographystyle{alpha}

\end{document}